\documentclass[12pt,reqno,final]{article}
\usepackage{chen}
\hypersetup{pdfpagemode=UseOutlines,colorlinks=false,pdfpagelayout=SinglePage,pdfstartview=FitH,bookmarksopen=true}
\begin{document}
\title{On the Existence of Global Bisections of Lie Groupoids}
\author{\textsc{Zhuo Chen} \\
{\small LMAM, Peking University \& CEMA, Central University of Finance and Economics,} \\
{\small
\href{mailto:chenzhuott@gmail.com}{\texttt{chenzhuott@gmail.com}}}
\and \textsc{Zhang-Ju Liu} \\
{\small  Department of Mathematics and LMAM, Peking University} \\
{\small \href{mailto:liuzj@pku.edu.cn}{\texttt{liuzj@pku.edu.cn}}}
\and \textsc{De-Shou Zhong}
\\
{\small Center of Mathematics, CYU}\\
{\small
\href{mailto:zhongdeshou@263.net}{\texttt{zhongdeshou@263.net}}} }
\date{}
\maketitle

\begin{abstract} We show that every source connected Lie groupoid always has global
bisections through any given point. This bisection can be chosen to
be the multiplication of some exponentials as close as possible to a
prescribed curve. The existence of bisections through more than one
prescribed points is also discussed. We give some interesting
applications of these results.



\paragraph{Funded by} CPSF(20060400017).

\end{abstract}

\section{Introduction}

The notion of groupoids generalizes that of both groups and the
Cartesian product of a set, namely the pair groupoid. It is
Ehresmann who first made the concept of groupoid (\cite{MR0095894,
MR0213410}) central to his vision of differential geometry. Lie
groupoids, originally called differentiable groupoids, were also
introduced by Ehresman (see several papers contained in
\cite{MR794193}) and theories of Lie groupoids, especially the
relationships with that of Lie algebroids defined by Pradines
(\cite{MR0214103,MR0216409}) have been investigated by many people
and much work has been done in this field. Readers can find the most
basic definitions and examples of (Lie) groupoids in the texts such
as Mackenzie's \cite{MR896907},
 his recent book \cite{MR2157566}, and \cite{Kumpera,MR1747916}. The importance of groupoid
theories were already shown in the studies on symplectic groupoids
and Poisson geometry, for example, illustrated by Weinstein
(\cite{MR866024,MR959095}), Coste (\cite{MR996653}), Dazord
(\cite{MR1104915}), Karas\"ev (\cite{MR1008479}), Zakrzewski
(\cite{MR1081011}) and many other authors.

By definition, a Lie groupoid is a groupoid where the set $Ob$ of
objects and the set $Mor$ of morphisms are both manifolds, the
source and target operations are submersions, and all the category
operations (source and target, composition, and identity-assigning
map) are smooth. In our humble opinion, a Lie groupoid can thus be
thought of as a ``many-object generalization'' of a Lie group, just
as a groupoid is a many-object generalization of a group.

On a group $G$, the left translations $L_g$: $h\mapsto gh$ form a
group which is isomorphic to $G$ itself. In the groupoid world, the
concepts of left translations are replaced by a family of elements
which is called a bisection of the groupoid \cite[I,1.4]{MR2157566}.
For a groupoid $\Groupoid$, we call any section $s\lon \Groupoid$ of
the $\alpha$-fibers for which $\beta\circ s$ is a diffeomorphism a
bisection of $\Groupoid$ (also known as admissible sections in
\cite{MR896907}). Bisections may be regarded as generalized elements
of the groupoids and similarly, the exponential maps take values in
the collections of all bisections. Although in most situations where
bisections are used, the question is local and the existence of
local bisection are of course correct, it has remained an unsettled
problem whether global bisection exist through an arbitrary point.
In this paper, we give an affirmative answer to this question.

We will prove that there exists a bisection through a given point of
a Lie groupoid which is source connected (Theorem \ref{Thm:Main}).
Moreover, we show that this bisection is of the form $\exp X_1 \exp
X_2 \cdots \exp X_k$, where $X_i$ are some sections of the
corresponding Lie algebroid and compact supported. Furthermore, we
prove that these $X_i$ can be chosen as close as possible to a
prescribed curve (Theorem \ref{Thm:MainStrong}).

We also prove that if the groupoid is transitive and the base space
is more than $1$-dimensional, then for any given points $g_1$,
$\cdots$, $g_n$, such that their sources and targets are subject to
the concordance condition (see (\ref{Eqt:concordance})), there
exists a bisection through all of them (Theorem
\ref{Thm:SeveralPoints}).

As stated by Mackenzie in \cite{MR2157566}: ``Groupoids possess many
of the features which give groups their power and importance, but
apply in situations which lack the symmetry which is characteristic
of group theory and its applications'', we shall apply our theorems
to varies cases of groupoids and obtain some interesting results.

This paper is organized as follows. Section \ref{Sec:LGLA} begins
with an account of basic concepts of Lie groupoids, its bisections,
Lie algebroids and exponentials (our conventions follow that of
\cite{MR2157566}). Section \ref{Sec:MainThmApp} gives the main
Theorem \ref{Thm:Main} (we also provide a stronger version of this
result, Theorem \ref{Thm:MainStrong}) and some applications are
added in. Section \ref{Sec:SeveralPoints} studies the existence
problem of bisections through more than one points. Section
\ref{Proof} is devoted to prove Theorem \ref{Thm:MainStrong} which
implies Theorem \ref{Thm:Main}, and a detailed proof is spilt into
several lemmas. To understand what is going on, the reader is
referred to the seven pictures illustrating the geometric images.

\section{Lie Groupoids and its Tangent Lie
algebroids}\label{Sec:LGLA}
\begin{defn}\cite{MR2157566}A groupoid consists of a set $\Groupoid$
and a subset $M\subset \Groupoid$, called respectively the groupoid
and the base, together with two maps $\alpha$ and $\beta$ from
$\Groupoid$ to $M$, called respectively the source and target, and a
partial multiplication $(g,h)\mapsto gh$ in $\Groupoid$ defined on
the set of composable pairs:
$$\Groupoid^{[2]}\defbe\set{(g,h)\in \Groupoid\times\Groupoid|\
\beta(g)=\alpha(h)},$$ all subject to the following conditions:
\begin{itemize}
\item[i)] $\alpha(gh)=\alpha(g)$ and $\beta(gh)=\beta(h)$ for all $(g,h)\in
\Groupoid^{[2]}$;
\item[ii)] $f(gh)=(fg)h$ for all $f,g,h\in \Groupoid$ such that
$\beta(f)=\alpha(g)$ and $\beta(g)=\alpha(h)$; \item[iii)]
$\alpha(x)=\beta(x)=x$ for all $x\in M$;
\item[iv)] $g\beta(g)=\alpha(g)g=g$ for all $g\in \Groupoid$;
\item[v)] each $g\in \Groupoid$ has a two-sided inverse $g\inverse$
such that $\alpha(g\inverse)=\beta(g)$, $\beta(g\inverse)=\alpha(g)$
and $g\inverse g=\beta(g)$, $gg\inverse=\alpha( g)$.
\end{itemize}
\end{defn}
A groupoid $\Groupoid$ on the base $M$, with respectively source and
target maps $\alpha$, $\beta$, will be denoted by $(\Groupoid
\alphabetaarrow{M})$, or, more briefly, $(\Groupoid, M)$.   We adopt
the convention that, whenever we write a multiplication $gh$, we are
assuming that is defined (see  Pic. 1))

For $x\in M$, its \textbf{orbit}, denoted by $O_x$, is the set
$\beta\circ\alpha\inverse(x)\subset M$.

Lie groupoid $(\Groupoid\alphabetaarrow{M})$ is a groupoid with
differential structures and the base space $M$ is an embedded
submanifold, the target and source maps are submersions and all the
operations are smooth. The orbit $O_x$ is also a submanifold of $M$.
We also notice that, for any $\alpha$-fiber $P=\alpha\inverse(x)$,
$\beta|_P: P\lon O_x$ is again a submersion.

We take the \emph{tangent Lie algebroid} $(Lie\Groupoid,\rho)$ of
$(\Groupoid,M)$ as
$$Lie\Groupoid\defbe \bigcup_{x\in M}T_x\alpha\inverse(x)=
\set{v\in T_x\Groupoid| x\in M, \alpha_*(v)=0}.$$ In turn, the
bracket of $\sections{Lie\Groupoid}$ is determined by the commutator
of left invariant vector fields and the anchor map is given by
$\rho=\beta_{\Groupoid*}|_M$.

{\unitlength=1mm
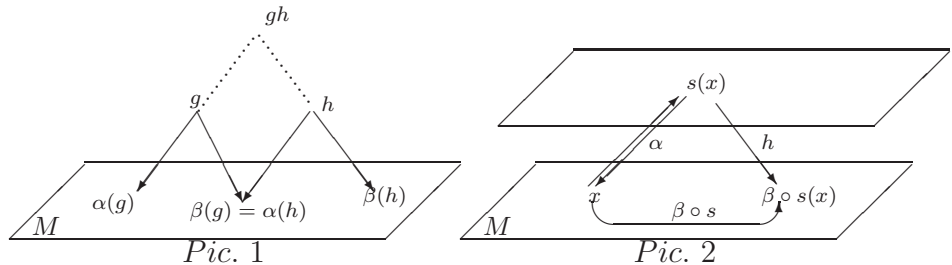
\begin{picture}(30,35)(0,0)

\put(15,3){\line(1,0){50.0}}\put(15,3){\line(1,1){10.0}}\put(18,3.5){\footnotesize
$M$} \put(25,13.2){\line(1,0){50.0}}\put(65,3){\line(1,1){10.0}}
 \multiput(47,29)(-0.6,-0.8){13}{\scriptsize$\cdot$}
 \put(40,20){\vector(-3,-4){8}}
  \put(26,7){\scriptsize$\alpha(g)$}
   \put(40,20){\vector(1,-2){6}}\put(39,21){\scriptsize$g$}
\put(39,6){\scriptsize$\beta(g)=\alpha(h)$}
   \multiput(47.8,29.4)(0.6,-0.8){12}{\scriptsize$\cdot$}
     \put(55.5,20){\vector(3,-4){8}}
  \put(62,8){\scriptsize$\beta(h)$}
\put(55,20){\vector(-3,-4){9}}\put(49,32){\scriptsize$gh$}\put(56.5,20){\scriptsize$h$}
\put(38,0){$Pic.~~~ 1$}

\put(75,3){\line(1,0){50.0}}\put(75,3){\line(1,1){10.0}}\put(78,3.5){\footnotesize
$M$} \put(85,13.2){\line(1,0){50.0}}\put(125,3){\line(1,1){10.0}}
\put(80,18){\line(1,0){50.0}}\put(80,18){\line(1,1){10.0}}
\put(90,28.2){\line(1,0){50.0}}\put(130,18){\line(1,1){10.0}}
\put(105,22){\vector(-1,-1){12}}\put(105,23){\scriptsize
$s(x)$}\put(100,15){\scriptsize$\alpha$}\put(92,8){\scriptsize$x$}\put(92,10){\vector(1,1){12}}
\put(103,6){\scriptsize$\beta\circ
s$}\put(109,21){\vector(3,-4){8}}\put(115,8){\scriptsize$\beta\circ
s(x)$}
\put(115,15){\scriptsize$h$}\put(105,8){\oval(25,6)[b]}\put(117.3,7){\vector(0,1){1}}
\put(98,0){$Pic.~~~ 2$}
\end{picture}}

\vskip 0.3in A bisection of a Lie groupoid $(\Groupoid,M)$ is a
smooth map $s: M\lon \Groupoid$ such that
\begin{itemize}
\item[1)] $\alpha\circ s=Id_M$;
\item[2)] $\beta\circ s$ is a diffeomorphism of $M$ (see Pic. 2).
\end{itemize}

The collection of all bisections of $\Groupoid$ is a group, in sense
of the following operations.
\begin{itemize}
\item[]\emph{Identity):}\\ The base $M$ serves as the identity, if regarded as a map
$M\lon \Groupoid$;
\item[]\emph{Multiplication):}\\ Let $s$ and $w$ be two bisections, their
multiplication $sw$ is again a bisection defined by (see Pic. 3)
$$sw(x)=s(x)w(\beta\circ s(x)),\quad\quad\forall x\in M;$$
\item[]\emph{Inversion):}\\ The inverse of $s$ is defined by (see Pic. 4)
$$s\inverse(x)=(s\circ(\beta\circ s)\inverse(x))\inverse,
\quad\quad\forall x\in M.$$
\end{itemize}

{\unitlength=1mm
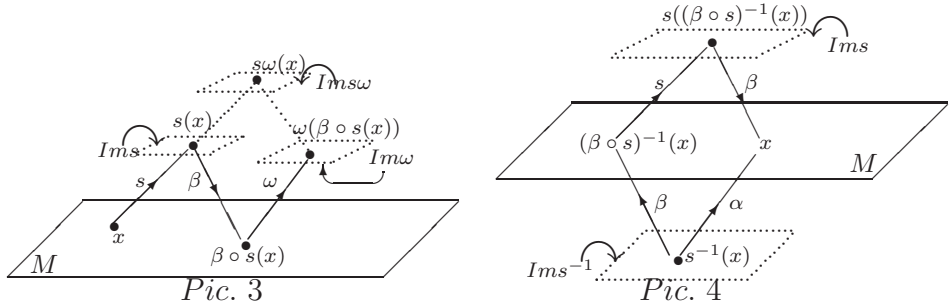
\begin{picture}(30,45)(0,0)

\put(15,3){\line(1,0){50.0}}\put(15,3){\line(1,1){10.0}}\put(18,3.5){\footnotesize
$M$} \put(25,13.2){\line(1,0){50.0}}\put(65,3){\line(1,1){10.0}}
\put(29,10){\line(1,1){10.0}}
\multiput(39.5,20)(0.5,0.5){17}{\scriptsize$\cdot$}
\put(29,10){\vector(1,1){6.0}}
\put(29,7.5){\scriptsize$x$}\put(28.5,9){\scriptsize$\bullet$}

\put(40,20){\line(1,-2){6}}\put(40,20){\vector(1,-2){3}}\put(37,23){\scriptsize$s(x)$}
\put(39,19.8){\scriptsize$\bullet$}
\put(30,20.5){\Large$\curvearrowright$}\put(27,19){\scriptsize$Ims$}

 \multiput(32,18.5)(0.8,0){13}{\scriptsize$\cdot$}
 \multiput(37,21)(0.8,0){12}{\scriptsize$\cdot$}
 \multiput(32,19)(0.8,0.4){6}{\scriptsize$\cdot$}
  \multiput(42,19)(0.8,0.4){6}{\scriptsize$\cdot$}
 \put(32,15){\scriptsize$s$}
\put(42,5){\scriptsize$\beta\circ
s(x)$}\put(39,15){\scriptsize$\beta$}
\put(46,6.5){\scriptsize$\bullet$}

 \multiput(48,17.5)(0.8,0){13}{\scriptsize$\cdot$}
  \multiput(53.5,20.5)(0.8,0){13}{\scriptsize$\cdot$}
  \multiput(48,17.6)(0.8,0.4){7}{\scriptsize$\cdot$}
   \multiput(58,17.6)(0.8,0.4){7}{\scriptsize$\cdot$}

   \put(49,15){\scriptsize$\omega$}
\multiput(48,28)(0.6,-0.8){12}{\scriptsize$\cdot$}
\put(61,17){\oval(8,3)[b]}\put(63,18){\scriptsize$Im\omega$}
\put(57,17){\vector(0,1){1}}

\multiput(40,27)(0.8,0){13}{\scriptsize$\cdot$}
  \multiput(45.5,29.5)(0.8,0){12}{\scriptsize$\cdot$}
   \multiput(40,27)(0.8,0.4){7}{\scriptsize$\cdot$}
    \multiput(50,27)(0.8,0.4){7}{\scriptsize$\cdot$}

\put(47,8){\line(3,4){8}}
\put(47,8){\vector(3,4){5}}\put(47.5,30.5){\scriptsize$s\omega(x)$}\put(47.5,28.5){\scriptsize$\bullet$}
\put(53,22){\scriptsize$\omega(\beta\circ
s(x))$}\put(54.5,18.5){\scriptsize$\bullet$}
\put(53,28.5){\Large$\curvearrowleft$}\put(56,28){\scriptsize$Ims\omega$}
\put(38,0){$Pic.~~~ 3$}

\multiput(93,2)(0.8,0){25}{\scriptsize$\cdot$}
 \multiput(93,2)(0.6,0.6){11}{\scriptsize$\cdot$}
   \multiput(99.6,8.5)(0.8,0){25}{\scriptsize$\cdot$}
   \multiput(113,2)(0.6,0.6){11}{\scriptsize$\cdot$}

\put(80,16){\line(1,0){50.0}}\put(80,16){\line(1,1){10.0}}\put(127,17){\footnotesize
$M$} \put(90,26.2){\line(1,0){50.0}}\put(130,16){\line(1,1){10.0}}
\put(96,21.5){\line(1,1){12}}\put(96,21.5){\vector(1,1){6}}
\put(101,12){\scriptsize$\beta$}\put(103,6){\line(-1,2){7}}\put(103,6){\vector(-1,2){4}}
\put(91,20){\scriptsize$(\beta\circ s)^{-1}(x)$}
\put(105,5){\scriptsize$s^{-1}(x)$}
\put(103.5,4.5){\scriptsize$\bullet$}
\put(104.5,6){\line(3,4){10}}\put(104.5,6){\vector(3,4){5}}
\put(115,21.5){\line(-1,2){6}}\put(112,28){\vector(1,-2){1}}\put(115,20){\scriptsize$x$}
\put(111,12){\scriptsize$\alpha$}\put(91,5){\Large$\curvearrowright$}\put(84,3){\scriptsize
$Ims^{-1}$}

\multiput(95,31.5)(0.8,0){25}{\scriptsize$\cdot$}
\multiput(95.5,32)(0.8,0.4){8}{\scriptsize$\cdot$}

 \multiput(101,35.2)(0.8,0){25}{\scriptsize$\cdot$}
  \multiput(115,32)(0.8,0.4){8}{\scriptsize$\cdot$}
\put(102,37){\scriptsize$s((\beta\circ s)^{-1}(x))$}
\put(108,33.5){\scriptsize$\bullet$}
\put(101,28){\scriptsize$s$}\put(113,28){\scriptsize$\beta$}

 \put(121,35){\Large$\curvearrowleft$}
\put(124,33){\scriptsize$Ims$}
 \put(99,0){$Pic.~~~ 4$}
\end{picture}}

\vskip 0.3in We recall the result of Kumpera and Spencer.
\begin{lem}\cite{MR0380908} Let $\Groupoid$ be any Lie groupoid
over base $M$, and let $\huaA=Lie \Groupoid$ be its Lie algebroid
with anchor $\rho$. For $X\in \sections{\huaA}$, let $\LeftMove{X}$
be the left invariant vector field corresponding to $X$, and recall
that $\rho(X)=\beta_*(\LeftMove{X})$ is its projected vector field
on M. Then $\LeftMove{X}$ is complete if and only if $\rho(X)$ is
complete. In fact, $\tilde{\phi}_t(g)$ is defined whenever
$\phi_t(\beta(g))$ is defined, where $\tilde{\phi}_t(g)$ and
${\phi}_t(g)$ are the flows generated by $\LeftMove{X}$ and
$\rho(X)$, respectively.
\end{lem}


As a corollary, let $X\in \sections{\huaA}$ be any section which has
a \textbf{compact support}, we know that $\LeftMove{X}$ is complete.
We denote $\exp(tX): M\lon \Groupoid$ ($t\in\Real$) the map
$$
\exp(tX)(x)\defbe \tilde{\phi}_t(x),\quad\quad\forall x\in M,
$$
called the exponential of $X$. One has
$$
\tilde{\phi}_t(g)=g\exp(tX)(\beta(g)),\quad\quad \forall
g\in\Groupoid.
$$ Furthermore,
$$\alpha\circ\exp(tX)=Id_M,\quad\quad\beta\circ\exp(tX)=\phi_t\,,$$
which shows that $\exp(tX)$ is a bisection of $\Groupoid$, for any
$t\in \Real$.

Another important fact is that
$$(\exp tX)\inverse = \exp(-tX).$$

In the special case that, the Lie groupoid degenerates to a group,
i.e., $M=pt$ is a point, the exponential map becomes the ordinary
exponential of Lie groups. We need the following basic results of
Lie group theories.

\begin{lem}\label{Lem:huaNandO}
Let $G$ be a Lie group and let $e\in G$ be its unit element. Then
there exists an open neighborhood $\huaN$ of $e$, and an open
neighborhood $\huaO\subset \LieG=T_eG$ near zero, such that
$$\exp: \huaO\lon \huaN$$
is a diffeomorphism \cite{MR514561,MR746308}.
\end{lem}

For a general Lie groupoid $(\Groupoid,M)$ and $x\in M$, by
$G_x\defbe \alpha\inverse(x)\cap\beta\inverse(x)$ we denote the
isotropic group at $x$. $G_x$ is a Lie group. By $\LieG_x\defbe
ker\rho_x\subset \huaA_x$ we denote the isotropic algebra at $x$,
which is a Lie algebra, in fact, the Lie algebra of $G_x$.

If $X\in \sections{\huaA}$ is a section with {compact support} and
$X_x\in \LieG_x$, then $\exp(tX)(x)$ is the usual exponential map
$\exp t(X_x)$, regarding $X_x\in \LieG_x$.

\begin{defn}A bisection of $(\Groupoid,M)$ is said to be   finitely
generated if it has the form $ \exp X_1 \exp X_2 \cdots \exp X_k$,
for some $X_1$, $\cdots$, $X_k\in\sections{\huaA}$ and every $X_i$
has a compact support. It is said to be finitely generated over
$\huaU$, an open set $\huaU\subset M$, if each support of $X_i$ is
contained in $\huaU$.
\end{defn}
\section{Main Theorems and their Applications}
\label{Sec:MainThmApp}

\begin{thm}\label{Thm:Main}
Let $(\Groupoid,M)$ be an $\alpha$-connected Lie groupoid. Then for
any $g\in\Groupoid$, there exists a finitely generated bisection $s$
of $\Groupoid$ through $g$, i.e.,  $s(\alpha(g))=g $.
\end{thm}

An immediate consequence of this theorem is the well known fact
that, for a connected Lie group $G$, every element $g$ can be
expressed into $$g=\exp X_1\exp X_2\cdots\exp X_k ,$$ for some
$X_i\in T_eG$.

In this paper we would like to prove a stronger version of Theorem
\ref{Thm:Main} stated as follows.
\begin{thm}\label{Thm:MainStrong}
Let $(\Groupoid,M)$ be Lie groupoid and $x\in M$. Let
$g\in\alpha\inverse(x)$ and let
$\widetilde{\huaU}\subset\alpha\inverse(x)$ be a connected open set
which contains both $x$ and $g$. Suppose that $\huaU$ is an open set
containing $\beta(\widetilde{\huaU})$,
then there exists a bisection $s$ of $\Groupoid$ through $g$, and
$s(y)=y$ for all $y\in \huaU^c$ ($=M-\huaU$). Moreover, $s$ is
finitely generated over $\huaU$.
\end{thm}

The proof of this theorem is given in the last section of this
paper. As an application of Theorem \ref{Thm:Main} as well as
\ref{Thm:MainStrong}, we have the following interesting conclusions.
We always assume that the reader is familiar with the various kinds
of groupoids mentioned below, especially their bisections.
\begin{thm}\label{Thm:Pair1}{\rm[Homogeneity of manifolds
\cite{MR1445290}.]}
Let $M$ be a connected smooth manifold. Then, for any two points
$x$, $y\in M$ and any connected open set $\huaU$ containing $x$ and
$y$, there exists a diffeomorphism $\Phi: M\lon M$, such that
$\Phi(x)=y$ and $\Phi(m)=m$ for all $m\in\huaU^c$.
\end{thm}

This theorem is of course implied by the following one.
\begin{thm}\label{Thm:Pair2} With the same assumptions as in the
previous one, for any two points $x$, $y\in M$, and for
an arbitrary open neighbor $\huaU$ containing $x$ and $y$, there
exist some smooth vector fields $X_1$, ..., $X_k$ compact supported
within $\huaU$, such that
$$\varphi_k^1\circ\cdots\circ\varphi_2^1\circ\varphi_1^1(x)=y,$$
where $\varphi_i^t$ is the flow of $X_i$.
\end{thm}

\noindent\textbf{Proof of Theorem \ref{Thm:Pair1} and
\ref{Thm:Pair2}.} We recommend \cite[Example 1.1.7,
1.4.3]{MR2157566} for background information on the pair groupoid
$M\times M$, for which a bisection is exactly a diffeomorphism of
$M$. The exponential of a vector field which has compact support is
its flow. And the multiplication of two bisections, namely
diffeomorphism, are exactly their compositions. The conclusion of
Theorem \ref{Thm:Pair2} is exactly the translation of Theorem
\ref{Thm:MainStrong} into the pair groupoid case. \qed

In what follows, we consider a vector bundle $E\stackrel{q}{\lon} M$
and we denote by $\Phi(E)$ the \textsl{linear frame groupoid} of
$E$, simply called the frame groupoid. Please refer to
\cite[III]{MR896907} and \cite[Example 1.1.12]{MR2157566}, where it
is denoted by $\Pi(E)$), which is the collection of all linear
isomorphisms from a fiber of $E$ to some generally different fiber
of $E$, i.e., an element in $\Phi(E)$ is a vector space isomorphisms
$\xi: E_x\lon E_y$ for $x,y\in M$. The bisections group of $\Phi(E)$
is in fact $Aut(E)$, the group of vector bundle automorphisms of
$E$. It is proved in that $\Phi(E)$ is also a Lie groupoid on $M$.
One may directly draw from Theorem \ref{Thm:Main} the following
result.

\begin{thm}\label{Thm:VectorBundle}
Let $(E\lon M)$ be a vector bundle over a connected smooth manifold
$M$. For any two points $x$, $y\in M$ and an isomorphism of vector
spaces $\phi: E_x\lon E_y$, there exists an automorphism $\Phi:
E\lon E$ of vector bundles such that $\Phi|_{E_x}=\phi$.
\end{thm}

Of course we can add some structures in $E$. If $(E, [~,~])$ is a
Lie algebra bundle, one has the Lie-algebra-bundle frame groupoid
$\Phi_{Aut}(E)$ (\cite[Example 1.7.12]{MR2157566}). If
$(E,\langle~,~\rangle)$ is a Riemannian vector bundle, one gets the
orthonormal frame groupoid $\Phi_{\huaO}(E)$ (\cite[Example
1.7.9]{MR2157566}). For these two examples, see also Corollary
3.6.11 in \cite{MR2157566}. We are then easy to draw the following
analogue conclusions.
\begin{thm}\label{Thm:LieAlgbraBundle}
Let $(E\lon M,[~,~])$ be a Lie algebra bundle over a connected
smooth manifold $M$. For any two points $x$, $y\in M$ and an
isomorphism of Lie algebras $\phi: E_x\lon E_y$, there exists an
automorphism $\Phi: E\lon E$ of Lie algebra bundles such that
$\Phi|_{E_x}=\phi$.
\end{thm}

\begin{thm}\label{Thm:RiemanBundle}
Let $(E\lon M, \langle~,~\rangle)$ be a Riemannian vector bundle
over a connected smooth manifold $M$. For any two points $x$, $y\in
M$ and an isomorphism of metric spaces $\phi: E_x\lon E_y$, there
exists an automorphism $\Phi: E\lon E$ of Riemannian vector bundles
such that $\Phi|_{E_x}=\phi$.
\end{thm}

Finally, we consider the action groupoid $M\sphericalangle G$
 coming from a right action of a connected Lie group
$G$ on a connected manifold $M$ (\cite[Example 1.1.9]{MR2157566},
see also \cite{MR896907}). Here $M\sphericalangle G=M\times G$ is a
Lie groupoid on $M$. Recall that a bisection of $M\sphericalangle G$
can be identified with a smooth $G$-valued function $s: M\lon G$
such that the map $$M\lon M,\quad m \mapsto ms(m),\quad \forall m\in
M$$ is a diffeomorphism of $M$. Such kinds of $s$ are called
\textbf{invertible} functions. So we have the following theorem.
\begin{thm}\label{Thm:ActionGroup} For any prescribed
$x\in M$, $g\in G$, one can find an invertible function $s: M\lon G$
satisfying $s(x)=g$.
\end{thm}

\section{Bisections through Points}\label{Sec:SeveralPoints}

Now we consider a more generalized problem: does there exists a
bisection through two (or more) given points of a groupoid? The
answer is also yes, under some topological conditions.

We recall that for a transitive Lie groupoid $(\Groupoid,M)$, the
map $\beta|_{\alpha\inverse(x)}: \alpha\inverse(x)\lon M$ is a
surjection as well as a submersion, for any $x\in M$.

The following one is the main theorem in this section. We always
assume that the natural number $n\geqslant 2$.
\begin{thm}\label{Thm:SeveralPoints}
Let $(\Groupoid,M)$ be a transitive and $\alpha$-connected Lie
groupoid and suppose that $dim M\geqslant 2$. For any different $n$
points $g_1$,$\cdots$,$g_n\in \Groupoid$, and let $\alpha(g_i)=x_i$,
$\beta(g_i)=y_i$, $i=1,\cdots,n$. Then there exists a bisection $s$
of $\Groupoid$ such that $s(x_1)=g_1$, $\cdots$, $s(x_n)=g_n$ if and
only if
\begin{equation}\label{Eqt:concordance}
x_i\neq x_j,\quad \mbox{and }\ y_i\neq y_j\,,\quad \forall i\neq j.
\end{equation}
Moreover, this bisection is  finitely generated.
\end{thm}
In what follows we devote to proving this theorem. We need some
preparations and let us introduce some concepts first. For $n$ pairs
of points
$$p_i=(x_i,y_i)\in M\times M,\quad i=1,\cdots,n,$$
 they are said to be \textbf{concordant} if they are subject to condition
(\ref{Eqt:concordance}). If they are concordant and there is a
proper arrangement of their indices such that they make a loop, we
say they consist a \textbf{chain}. That is, for some permutation of
the indices $\check{p}_i=p_{\theta(i)}=(\check{x}_i,\check{y}_i)
=(x_{\theta{(i)}},y_{\theta{(i)}})$, where $\theta\in \huaS(n)$ (the
permutation group), one has
$$
\check{y}_1=\check{x}_2,\quad
\check{y}_2=\check{x}_3,\quad\cdots,\quad \check{y}_n=\check{x}_1.
$$
We shall write
$$
\check{p}_1 \curvearrowright \check{p}_2 \curvearrowright \cdots
\curvearrowright \check{p}_n
$$
to denote such a chain (see Pic. 5).

\vskip 0.3in
\begin{picture}(30,35)(-40,0)
\put(0,0){\vector(0,1){18}} \put(-35,0){$\check{x}_1=\check{y}_n$}
\put(0,18){\vector(1,1){16}} \put(-30,20){$\check{x}_2=\check{y}_1$}
\put(16,34){\scriptsize$\cdots\quad\cdots\quad$}
\put(50,34){\vector(1,-1){30}}
\put(85,10){$\check{x}_{n-1}=\check{y}_{n-2}$}
\put(80,4){\vector(-1,-1){30}} \put(50,-26){\vector(-2,1){50}}
\put(55,-30){$\check{x}_{n}=\check{y}_{n-1}$} \put(-20,-30){$Pic.~~~
5$}

\put(200,-20){\vector(1,1){28}} \put(229,8){\vector(1,0){38}}
\put(222,30){\vector(1,-1){48}} \put(270,-18){\vector(1,0){28}}
\put(269,8){\vector(1,1){28}} \put(210,-30){$Pic.~~~ 6$}

\end{picture}
\vskip 0.6in

\begin{defn} Let $p_i=(x_i,y_i)\in M\times M$, $i=1,\cdots,n$
be some pairs of points in $M$. $p_1$, $\cdots$, $p_n$ are said to
be \textbf{independent}, if they are concordant and there is not any
subset of $\set{p_i}$ that can consist a chain.
\end{defn}
\begin{rmk} { We recall
the pair groupoid $M\times M$ of pairs of points $p=(x,y)$. We say
an elements $p\in M\times M$ can be expressed by some $p_1$,
$\cdots$, $p_m$ $\in M\times M$, if there are some $m_1$,$\cdots$,
$m_k\in \set{1,\cdots,m}$, such that
$$p=q_{m_1}q_{m_2}\cdots q_{m_k},$$
where each $q_{m_l}$ is $p_{m_l}$ or the inverse $p_{m_l}\inverse$.
One can prove that, for concordant $n$ elements $p_1$, $\cdots$,
$p_n$, they are independent if and only if each $p_i=(x_i,y_i)$ is
not possible to be expressed by
${p_1,p_2,\cdots,\widehat{p_i},\cdots,p_n}$. This is the reason that
we use ``independent''.}
\end{rmk}
For example, the five pairs in Pic.6 are independent.

\begin{lem}\label{Lem:InddtWellOder}
Let $p_i=(x_i,y_i)\in M\times M$, $i=1,\cdots,n$ be some pairs of
points in $M$. Then $p_1$, $\cdots$, $p_n$ are independent if and
only if there exists some $\sigma\in \huaS(n)$, such that for
$\bar{p}_1=(\bar{x}_1,\bar{y}_1)=p_{\sigma(1)}=(x_{\sigma(1)},y_{\sigma(1)})$,
$\cdots$,
$\bar{p}_n=(\bar{x}_n,\bar{y}_n)=p_{\sigma(n)}=(x_{\sigma(n)},y_{\sigma(n)})$,
one has
\begin{itemize}
\item[] $\bar{x}_2\neq \bar{y}_1$;
\item[] $\bar{x}_3\neq \bar{y}_1$, $\bar{x}_3\neq \bar{y}_2$;
\item[] $\cdots\cdots$,
\item[] $\bar{x}_k\neq \bar{y}_l$, for $l {=1,\cdots,k-1}$;
\item[] $\cdots\cdots$,
\item[] $\bar{x}_n\neq \bar{y}_l$, for $l
{=1,\cdots,n-1}$.
\end{itemize}
(In this case and for convenience, we will say $\bar{p}_{i}$ are
well-ordered.)
\end{lem}

\pf ``$\Leftarrow$'': We adopt a negative approach. If there are
some elements $\bar{p}_{m_1}$, $\cdots$, $\bar{p}_{m_k}$ consist a
chain
$$
\bar{p}_{m_1}\curvearrowright \bar{p}_{m_2} \curvearrowright\cdots
\curvearrowright \bar{p}_{m_k},
$$
then from $\bar{y}_{m_1}=\bar{x}_{m_2}$, we conclude $m_1> m_2$.
Similarly, from $\bar{y}_{m_j}=\bar{x}_{m_{j+1}}$, we conclude
$m_{j}>m_{j+1}$ and finally one has $m_1>m_2>\cdots>m_k$. On the
other hand, $\bar{x}_{m_1}=\bar{y}_{m_k}$ implies $m_k>m_1$:
contradiction!





``$\Rightarrow$'': We give an inductive proof. For $n=2$, if two
elements $p_1=(x_1,y_1)$ and $p_2=(x_2,y_2)$ are independent, then
$p_1\neq p_2\inverse$, i.e., either $x_2\neq y_1$ or $x_1\neq y_2$
holds. So one can always find $\sigma\in\huaS(2)$, which is either
$Id$ or the flip $(1,2)$.

Suppose that for $n\geqslant 2$, the lemma holds. For any
independent $n+1$ elements $p_1$, $\cdots$, $p_{n+1}$, we claim that
there exist some $k\in \set{1,\cdots, n+1}$, such that
$$x_k\neq y_j,\quad \forall j\in \set{1,\cdots,\widehat{k},\cdots,n+1}.$$
In fact, if it is not true, then for each $x_i$, one find some
$\psi(i)\in \set{1,\cdots,\widehat{i},\cdots,n+1}$ and
$x_i=y_{\psi(i)}$. Obviously $\psi(i)$ is unique. It is also easy to
see that $\psi$ is a permutation of $n+1$ numbers.  Thus, if we find
the smallest number $m\geqslant 1$ such that $\psi^{m+1}(1)=1$, then
we find a chain:
$$
 (x_{\psi^m(1)},y_{\psi^m(1)})\curvearrowright
 (x_{\psi^{m-1}(1)},y_{\psi^{m-1}(1)})\curvearrowright\cdots
 \curvearrowright(x_1,y_1),
$$
which contradicts with the assumption that $p_1$, $\cdots$,
$p_{n+1}$ are independent.

By this claim we pick $\sigma_0=(k,n+1)\in\huaS(n+1)$, which is the
flip of $k$ and $n+1$ and now for
$$\tilde{p}_{i}=(\tilde{x}_i,\tilde{y}_i)=
p_{\sigma_0(i)}=(x_{\sigma_0(i)},y_{\sigma_0(i)}), i=1,\cdots,
n+1,$$ one has
\begin{equation}\label{temp:2}
\tilde{x}_{n+1}=x_{k}\neq \tilde{y}_{j}=y_{\sigma_0(j)},\quad
\forall j=1,\cdots, n. \end{equation}

By the inductive assumption, we are able to find
$\sigma_1\in\huaS(n)$, such that
$$
\bar{p}_i=\tilde{p}_{\sigma_1(i)}, \quad i=1,\cdots,n,
$$
satisfies the well-ordered condition. Of course if we write
$\bar{p}_{n+1}=\tilde{p}_{n+1}$, then (\ref{temp:2}) shows that
$\bar{p}_1$, $\cdots$, $\bar{p}_{n+1}$ are also well-ordered. \qed


\begin{lem}\label{Lem:Seperate}
Let $(\Groupoid,M)$ be a transitive and $\alpha$-connected Lie
groupoid. Suppose that $dim M\geqslant 2$ and let $x_1$, $\cdots$,
$x_k$ be some points of $M$. Then for any $g\in \Groupoid$, if
$x=\alpha(g)$ and $y=\beta(g)$ are not contained in the set
$\cup_{i=1}^k\set{x_i}$, then there exists a bisection $s$ such that
$s(x)=g$ and $s(x_i)=x_i$, for all $i=1$, $\cdots$, $k$.
\end{lem}
\pf Consider the open set
$\widetilde{\huaU}=(\beta\inverse(x_1)\cup\cdots\cup\beta\inverse(x_k))^c\cap
\alpha\inverse(x)$. Obviously $x$ and $g$ are contained in
$\widetilde{\huaU}$. Since $\beta|_{\alpha\inverse(x)}:
\alpha\inverse(x)\lon M$ is a surjective submersion, we know that
$$dim \alpha\inverse(x)-dim(\beta\inverse(x_i)\cap
\alpha\inverse(x))=dim M\geqslant 2.$$ And hence by
$\alpha\inverse(x)$ being connected, $\widetilde{\huaU}$ is also
connected. Let
$\huaU=\beta(\widetilde{\huaU})=(\cup_{i=1}^k\set{x_i})^c$ be the
corresponding open set in $M$. So Theorem \ref{Thm:MainStrong}
claims that there exists a bisection $s$ with $s(x)=g$ and
$s(x_i)=x_i$, $i=1$, $\cdots$, $k$. \qed

\begin{prop}\label{Pro:Independent}
Let $(\Groupoid,M)$ be a transitive and $\alpha$-connected Lie
groupoid and suppose that $dim M\geqslant 2$. For $n$ points
$g_1$,$\cdots$,$g_n\in \Groupoid$, and let $\alpha(g_i)=x_i$,
$\beta(g_i)=y_i$, $i=1,\cdots,n$, if
$$p_i=(x_i, y_i),\quad i=1,\cdots,n,$$
are independent, then there exists a bisection $s$ of $\Groupoid$
such that $s(x_1)=g_1$, $\cdots$, $s(x_n)=g_n$.
\end{prop}
\pf By Lemma \ref{Lem:InddtWellOder}, it suffices to assume that
those $p_i$ are already well-ordered. I.e., for each $i$,
$$
x_i \mbox{ and } y_i \notin \bigcup_{j=1}^{i-1}\set{y_j} \cup
\bigcup_{k=i+1}^n\set{x_k}.
$$
Thus, Lemma \ref{Lem:Seperate} tells us that there exists a
bisection $s_i$ such that
$$
s_i(x_i)=g_i,\quad s_i(y_j)=y_j,\quad\forall j<i,\ \mbox{ and }\
s_i(x_k)=x_k,\quad \forall k>i.
$$
Now, let $w_1=s_1$, $w_{2}=s_1 s_2$, $\cdots$, $w_{n}=s_1s_2\cdots
s_n$ and we claim $s=w_n$ is the bisection we are looking for. We
show this fact inductively. Of course we have
$$w_1(x_1)=g_1,\quad w_1(x_j)=x_j, \ j= 2,\cdots,n.$$ Suppose
that
$$w_k(x_j)=g_j,\ \forall j\leqslant k,\quad
w_k(x_l)=x_l,\ \forall l\geqslant{k+1}
$$
is already proved, then for $j<k+1$,
\begin{eqnarray*}
&&w_{k+1}(x_j) = (w_k s_{k+1})(x_j)\\
&=& w_k(x_j) s_{k+1}(\beta\circ w_k (x_j))= g_j s_{k+1}(y_j)=g_j
y_j=g_j.
\end{eqnarray*}
For $k+1$,
\begin{eqnarray*}
&&w_{k+1}(x_{k+1}) = (w_k s_{k+1})(x_{k+1})\\
&=& w_k(x_{k+1}) s_{k+1}(\beta\circ w_k (x_{k+1}))= x_{k+1}
s_{k+1}(x_{k+1})=g_{k+1}.
\end{eqnarray*}
And for $l>k+1$,
\begin{eqnarray*}
&&w_{k+1}(x_l) = (w_k s_{k+1})(x_l)\\
&=& w_k(x_l) s_{k+1}(\beta\circ w_k (x_l))= x_l s_{k+1}(x_l)=x_l.
\end{eqnarray*}
This completes the proof.\qed

\noindent\textbf{Proof of Theorem \ref{Thm:SeveralPoints}.}
``$\Rightarrow$'': If there exists a bisection $s$ of $\Groupoid$
such that $s(x_1)=g_1$, $\cdots$, $s(x_n)=g_n$ and these $g_1$,
$\cdots$, $g_n$ are different points in $\Groupoid$, then
$$g_i= s(x_i),\quad g_i\inverse=s\inverse(y_i).$$
Since both $s$ and $s\inverse$ are both embeddings of $M$ into
$\Groupoid$ and $g_i\neq g_j$ ($i\neq j$), $x_i\neq x_j$, $y_i\neq
y_j$, this shows that these $p_i=(x_i,y_i)$ are concordant $n$ pairs
of points in $M$.

``$\Leftarrow$'': Suppose that these $p_i$ are concordant. We find
all subsets of $\set{p_1,\cdots,p_n}$ which consist chains. It is
easy to see that any two such subsets are disjoint. Let $k$ be the
number of these chains.

We induct with $k$, which is obviously less than $\half {(n+1)}$. If
$k=0$, i.e., these $p_i$ are independent, Proposition
\ref{Pro:Independent} already assures the existence of such a
bisection through all $g_i$. We assume that for $k\leqslant m$ the
conclusion is right. Now suppose that one has all the chains $S_1$,
$\cdots$, $S_{m+1}$, it suffices to assume $S_{m+1}$ is the chain
$$
p_1\curvearrowright p_2 \curvearrowright \cdots \curvearrowright
p_r,
$$
where $r=\sharp S_{m+1}\leqslant n$. Choose an arbitrary $y_0\in M$
such that
$$y_0\notin
\bigcup_{j=1}^{n}\set{x_j} \cup \bigcup_{j=1}^n\set{y_j}.
$$
Since $\Groupoid$ is transitive, one can find some $h_0\in
\Groupoid$ such that $\alpha(h_0)=y_0$, $\beta(h_0)=x_1$. Let $g_0=
g_r h_0\inverse$. Hence $\alpha(g_0)=x_r$, $\beta(g_0)=y_0$. It is
easy to see that $g_0$, $g_1$, $g_2$, $\cdots$, $\widehat{g_r}$,
$\cdots$, $g_{n}$ are different points in $\Groupoid$, and the
corresponding
$$
p_0=(x_r,y_0),
p_1=(x_1,y_1),\cdots,\widehat{p_{r}},\cdots,p_{n}=(x_{n},y_{n})\in
M\times M
$$
has only $m$ chains $S_1$, $\cdots$, $S_{m}$. And we are able to
determine a bisection $\tilde{s}$ such that
$$
\tilde{s}(x_r)=g_0, \quad \tilde{s}(x_i)=g_i,\quad \forall i\in
\set{1,\cdots,\widehat{r},\cdots,n}.
$$
By the mean time, Lemma \ref{Lem:Seperate} gives a bisection
$\breve{s}$ such that
$$\breve{s}(y_0)=h_0, \quad \breve{s}(y_i)=y_i,
\quad \forall i\in \set{1,\cdots,\widehat{r},\cdots,n}.
$$
Then it is a direct check that the bisection $\tilde{s}\breve{s}$
satisfies
\begin{eqnarray*}
&& \tilde{s}\breve{s} (x_i)\quad (i\neq r)\\
&=& \tilde{s}(x_i)\breve{s}(\beta\circ \tilde{s}(x_i)) = g_i
y_i=g_i,
\end{eqnarray*}
and
\begin{eqnarray*}
&& \tilde{s}\breve{s} (x_r) = \tilde{s}(x_r)\breve{s}(\beta\circ
\tilde{s}(x_r)) = g_0 \breve{s}(x_0)=g_0 h_0=g_r.
\end{eqnarray*}
This shows that $\tilde{s}\breve{s}$ is just what we need. The
preceding process of construction of such a bisection also shows
that $\tilde{s}\breve{s}$ is finitely generated. This completes the
proof.
 \qed

\begin{rmk}\rm From the proof, one is able to see that, in Theorem
\ref{Thm:SeveralPoints}, the condition ``$\Groupoid$ is transitive
and $dim M\geqslant 2$'' can be replaced by ``each orbit of the base
space $M$ is more than $1$-dimensional''.
\end{rmk}

In what follows we present some applications of Theorem
\ref{Thm:SeveralPoints} in certain kind of groupoids. They are
respectively similar to Theorem \ref{Thm:Pair1},
\ref{Thm:VectorBundle}, \ref{Thm:LieAlgbraBundle},
\ref{Thm:RiemanBundle} and \ref{Thm:ActionGroup} but concerning
several pairs of points, and the proofs are omitted.

\begin{thm}\label{Thm:Pair3}
Let $M$ be a connected smooth manifold. Then, for any prescribed
points $x_i$, $y_i\in M$ such that $(x_i,y_i)$, $i=1,\cdots,n$ are
concordant, there exists a diffeomorphism $\Phi: M\lon M$, such that
$\Phi(x_i)=y_i$,  $i=1,\cdots,n$.
\end{thm}

\begin{thm}\label{Thm:VectorBundle2}
Let $(E\lon M)$ be a vector (resp. Lie algebra, Riemannian) bundle
over a connected smooth manifold $M$. For any prescribed points
$x_i$, $y_i\in M$ and isomorphisms of vector (resp. Lie algebra,
Riemannian) spaces $\phi_i: E_{x_i}\lon E_{y_i}$, such that
$(x_i,y_i)$, $i=1,\cdots,n$ are concordant, there exists an
automorphism $\Phi: E\lon E$ of vector (resp. Lie algebra,
Riemannian) bundles such that $\Phi|_{E_{x_i}}=\phi_i$.
\end{thm}

\begin{thm}\label{Thm:ActionGroup2}
Let $G$ be a connected Lie group which acts on a manifold $M$ and
suppose that $dim M\geqslant 2$. If the action is transitive, then
for any prescribed $x_i\in M$, $g_i\in G$ such that $(x_i,x_ig_i)$,
$i=1,\cdots,n$ are concordant, then one can find an invertible
function $s: M\lon G$ satisfying $s(x_i)=g_i$, $i=1,\cdots,n$.
\end{thm}

\section{Proof of Theorem \ref{Thm:MainStrong}}
\label{Proof}

We split the proof of Theorem \ref{Thm:MainStrong} into several
steps. In this section, we fix a Lie groupoid $(\Groupoid,M)$ which
is $\alpha$-connected. We also assume that the base space $M$ is
connected.  Let $\huaA=Lie\Groupoid$.

\begin{lem}\label{Lem:WellCurve}
 Let $x\in M$, $g\in \alpha\inverse(x)$ and $c: [0,1]\lon
 \alpha\inverse(x)$ be a smooth curve such that $c(0)=x$, $c(1)=g$.
 If the base curve $\bar{c}=\beta\circ c: [0,1]\lon M$ is an
 injection and suppose that $\bar{c}$ is contained in some open set $\huaU\subset M$,
 then there exists some $X\in \sections{\huaA}$ with compact support in
 $\huaU$, such that $$\exp tX(x)=c(t),\quad \forall t\in [0,1].$$
 In particular, $\exp X$ is a bisection through $g$ and $\exp X|_{\huaU^c}$ is the identity map.
\end{lem}
\pf For each $t\in [0,1]$, we define $X_t\in \huaA_{\bar{c}(t)}$ to
be
$$X_t=l_{c(t)\inverse*} c'(t). $$
Since $\bar{c}$ does not intersect with itself, one is able to
extend this $t$-function into a well defined section
$X\in\sections{\huaA}$ which is compact supported in $\huaU$. It is
clearly that $\exp tX(x)=c(t)$.\qed

\begin{lem}\label{Lem:Local}
For each $x\in M$ and each open neighborhood $\huaU$ near $x$, there
exists an open set $\huaW\subset \alpha\inverse(x)$ containing $x$,
such that for each element $g\in \huaW$, there exists a bisection
$s: M\lon \Groupoid$ through $g$ and $s(y)=y$ for all $y\in
\huaU^c$. Moreover, $s$ has the form $s=\exp X$, for some
$X\in\sections{\huaA}$ and $X$ has a compact support contained in
$\huaU$.
\end{lem}
\pf The target map $\beta: \Groupoid\lon M$ is a submersion. For the
$\alpha$-fiber $P_x=\alpha\inverse(x)$, $\beta|_{P_x}: P_x\lon O_x$
is also a submersion. Let $dim O_x=m$, $dim P_x=m+n$. Hence, there
exist two local coordinate systems
$(\huaS;x_1,\cdots,x_m,y_1,\cdots,y_n)$ of $P_x$ near $x$ and
$(\huaT;x_1,\cdots,x_m)$ of $O_x$ near $x$ and they are subject to
the following requirements:
\begin{itemize}
\item[1)]  $\huaS\cong \Real^{m+n}$, $x$ is the origin point;
\item[2)]  $\huaT\cong \Real^{m}$, $x$ is also the origin point;
\item[3)] $\beta|_{\huaS}$ is given canonically by
$$\beta: (x_1,\cdots,x_m,y_1,\cdots,y_n)\mapsto
(x_1,\cdots,x_m).$$
\item[4)] $\huaT$ is contained in $\huaU$
and $\overline{\huaT}\subset \huaU$ is compact.
\end{itemize}

So, for two different points $g$, $h\in\huaS$, if $\beta(g)\neq
\beta(h)$, one is able to find a curve $c$ connecting $g$ and $h$.
Moreover, this curve can be chosen so that it  lies entirely in
$\huaS$, such that $\bar{c}=\beta\circ c$, which lies in $\huaT$, is
a curve without self-intersections.


Let $G_x=\alpha\inverse(x)\cap\beta\inverse(x)$ be the isotropic
group at $x$. Find the two open sets $\huaN$ and $\huaO$ as claimed
by Lemma \ref{Lem:huaNandO}. Of course we can assume that they are
both simply connected.

Write $G_0=G_x\cap \huaS$, which is a closed subset of $\huaS$.
Write $\huaN_0=\huaN\cap \huaS$, which is an open set of $G_0$,
where $G_0$ has the relative topology coming from $\huaS$.
Therefore, for each point $h\in \huaN_0$, one is able to find an
open ball $B(h)$ of $\huaS$ with $h$ at the center, such that
$$B(h)\cap \huaN_0\subset \huaN_0\,.$$

Now, let $$\huaW\defbe \cup_{h\in \huaN_0} B(h)\subset \huaS.$$ We
claim that this $\huaW$ is just what we need. In fact, for each
$g\in \huaW$, there are possibly two cases:

\emph{Case 1)} $\beta(g)=x$, i.e., $g\in G_x\cap \huaW=\huaN_0$.
Then by Lemma \ref{Lem:huaNandO}, we find an $X_x\in \huaO\subset
\huaA_x$, such that $\exp(X_x)=g$. Then extend $X_x$ arbitrarily to
a section $X\in \sections{\huaA}$ with compact support contained in
$\huaU$ and we get a bisection $\exp X$ through $g$.

\emph{Case 2)} $\beta(g)\neq x$. In this case, one can of course
find a smooth curve $c:[0,1]\lon \huaW$, such that $c(0)=x$,
$c(1)=g$, and more importantly, the curve $\bar{c}=\beta\circ c$
which connects $x$ with $\beta(g)$, is a curve who does not
intersect with itself. Then by Lemma \ref{Lem:WellCurve}, we also
obtain an $X\in \sections{\huaA}$ with compact support contained in
$\huaU$ and we get a bisection $\exp X$ through $g$.

Both of the two kinds of $X$ we constructed vanish outside of
$\huaU$. Hence the bisection $\exp X$ maps every $y\in\huaU^c$ to
itself. \qed

With these preparations, we are able to prove the strong version
main theorem.

\noindent{\textbf{Proof of Theorem \ref{Thm:MainStrong}}}. We choose
a smooth curve $c:[0,1]\lon \alpha\inverse(x)$ which lies in
$\widetilde{\huaU}$ such that $c(0)=\alpha(g)$, $c(1)=g$, and hence
the open set $\huaU$ covers the base curve $\bar{c}=\beta\circ c$
which connects $\alpha(g)$ and $\beta(g)$.

For each $t\in [0,1]$, consider the point $\bar{c}(t)\in M$. Set
$\huaU=M$, Lemma \ref{Lem:Local} says that there is neighborhood
$\huaW_t\subset \alpha\inverse(\bar{c}(t))$ near $\bar{c}(t)$, such
that for each $h\in\huaW_t$, there is a bisection through $h$.

Now, these $l_{c(t)}\huaW_t$ become an open coverage of the curve
$c$. Since the interval $[0,1]$ is compact, we find the following
data (see Pic. 7):
\begin{itemize}
\item[1)] a partition of $[0,1]$:
$$0=t_0<t_1<t_2<\cdots<t_{k-1}<t_k=1;$$
\item[2)] finitely some open sets $\huaW_i=\huaW_{t_i}$,
$i=0,1,\cdots,k$, such that $\bar{c}(t_i)\in \huaW_i$; and hence the
collection of $l_{c(t_i)}\huaW_i$ ($i=0,1,\cdots,k$) finitely covers
the curve $c$;
\item[3)] some points
$$a_1\in (t_0,t_1),\quad a_2\in (t_1,t_2),\quad\cdots,\quad a_k\in
(t_{k-1},t_k),$$ such that
\begin{eqnarray*}&&c(a_1)\in \huaW_0\cap l_{c(t_1)}\huaW_1,\quad
c(a_2)\in l_{c(t_1)}\huaW_1\cap
l_{c(t_2)}\huaW_2,\\
&&\qquad\qquad\qquad\qquad\cdots,\quad c(a_{k})\in
l_{c(t_{k-1})}\huaW_{k-1}\cap l_{c(t_k)}\huaW_k\,.
\end{eqnarray*}
\end{itemize}

{\unitlength=1mm
\begin{picture}(0,35)(0,0)
\put(75,6){\line(0,1){10}}\put(15,16){\oval(120,16)[tr]}
\put(40,24){\circle{45}}\put(39.5,23){$\bullet$}\put(35,33){$l_{c(t_i)}W_i$}\put(37,21){\scriptsize$c(t_i)$}
\put(44,23){$\bullet$}\put(44.5,14){\vector(0,1){9.0}}\put(43,12){\scriptsize$c(a_i)$}
\put(50,24){\circle{45}}\put(49.5,23){$\bullet$}\put(57,28){$l_{c(t_{i-1})}W_{i-1}$}\put(48,21){\scriptsize$c(t_{i-1})$}
\put(15,24){\circle{35}}\put(14.5,23){$\bullet$}
\put(13,33){$l_gW_k$}\put(11,21){\scriptsize$g=c(1)$}

\put(75,6){\circle{35}}\put(74.5,5){$\bullet$}\put(81,10){$W_0$}\put(65,3){\scriptsize$\alpha(g)=c(0)=\bar{c}(0)$}

\put(75,6){\line(1,0){65}}\put(105,6){\circle{40}}\put(104.5,5){$\bullet$}
\put(105.2,6){\line(0,1){4}}\put(104.7,9){\scriptsize$\bullet$}\put(100,18){\scriptsize$c(t_{i-1})^{-1}c(a_i)$}
\put(109,18){\vector(-1,-2){3.5}} \put(95,7){\scriptsize$s_i$}
\put(102,3){\scriptsize$\bar{c}(t_{i-1})$} \put(91,1){$W_{i-1}$}
 \multiput(97,6)(0.6,0.6){4}{\scriptsize$\cdot$}
\multiput(99,8)(0.6,0.3){4}{\scriptsize$\cdot$}
\multiput(102,9)(0.6,0){12}{\scriptsize$\cdot$}
\multiput(110,8.5)(0.6,-0.6){4}{\scriptsize$\cdot$}
\put(115,6){\circle{40}}\put(114.5,5){$\bullet$}
\put(114,7){\scriptsize$\bar{c}(t_i)$}
\put(115.2,5.5){\line(0,-1){4}} \put(114.7,1){\scriptsize$\bullet$}
\multiput(112,1)(0.6,0){12}{\scriptsize$\cdot$}
\multiput(112,1)(-0.5,0.4){7}{\scriptsize$\cdot$}
\multiput(119,1)(0.5,0.4){7}{\scriptsize$\cdot$}
\put(123,3){\scriptsize$\omega_i$}
\put(110,-1){\scriptsize$c(t_i)^{-1}c(a_i)$} \put(121,10){$W_i$}
\put(140,6){\circle{40}}\put(139.5,5){$\bullet$}
\put(135,3){\scriptsize$\beta(g)=\bar{c}(1)$} \put(129,1){$W_k$}

\put(60,-5){\line(1,0){80.0}}\put(60,-5){\line(1,2){10.0}}
\put(70,15){\line(1,0){80.0}}\put(140,-5){\line(1,2){10.0}}
\put(65,-2.5){\footnotesize $M$}

\put(30,0){$Pic.~~~ 7$}
\end{picture}}

\vskip 0.3in By the last condition, we know that there is a
bisection $s_1$ through $c(a_1)$. And since
$l_{c(t_1)\inverse}c(a_1)\in \huaW_1$, there exists a bisection
$w_1$ through $l_{c(t_1)\inverse}c(a_1)$. Similarly, since
$l_{c(t_{i-1})\inverse}c(a_i)\in \huaW_{i-1}$, there is a bisection
$s_i$ through $l_{c(t_{i-1})\inverse}c(a_i)$. And since
$l_{c(t_{i})\inverse}c(a_i)\in \huaW_{i}$, there exists a bisection
$w_i$ through $l_{c(t_{i})\inverse}c(a_i)$. $\cdots\cdots$ We can
require these bisections $s_1$, $\cdots$, $s_{k}$, and $w_1$,
$\cdots$, $w_k$ are the identity maps on $\huaU^c$ and they are all
of the forms $s_i=\exp X_i$, $w_i=\exp Y_i$, for some $X_i$, $Y_i\in
\sections{\huaA}$ with compact supports in $\huaU$.

We now show that the section $s_1w_1\inverse s_2w_2\inverse\cdots
s_{k}w_k\inverse$ is a bisection through $g$.  We notice the
following inductive formulas, hold for all $i=1,\cdots,k$:
\begin{itemize}
\item[1)] $\beta\circ s_{i}(\bar{c}(t_{i-1}))=\bar{c}(a_i)$;
\item[2)] $\beta\circ w_{i}(\bar{c}(t_{i}))=\bar{c}(a_i)$.
\end{itemize}
Using these, one is able to get
\begin{itemize}
\item[3)] $w_i\inverse(\bar{c}(a_i))=w_i(\bar{c}{(t_i)})\inverse=
c(a_i)\inverse c(t_i)$;
\item[4)]$(s_iw_i\inverse)(\bar{c}(t_{i-1}))=l_{c(t_{i-1})\inverse}c(t_i) =
{c(t_{i-1})\inverse}c(t_i)$.
\end{itemize}
And hence we  obtain
\begin{eqnarray*}
&&s_1w_1\inverse s_2w_2\inverse\cdots
s_{k}w_k\inverse(x)\\
&=& s_1w_1\inverse(\bar{c}(t_0)) s_2w_2\inverse(\bar{c}(t_1)) \cdots
s_iw_i\inverse(\bar{c}(t_{i}-1)) \cdots
s_kw_k\inverse(\bar{c}(t_{k}-1))
\\
&=& c(t_1) c(t_1)\inverse c(t_2) \cdots {c(t_{i-1})\inverse}c(t_i)
\cdot {c(t_{k-1})\inverse}c(t_k)\\
&=& c(t_k)=g.
\end{eqnarray*}

It is also easy to check that  $s(y)=y$ for all $y\in \huaU^c$.
Since $w_k\inverse=\exp(-Y_i)$, we know that $s_1w_1\inverse
s_2w_2\inverse\cdots s_{k}w_k\inverse$ is just the bisection we
need. This completes the proof. \qed
\vskip0.4in

\noindent\textbf{Acknowledgements}

We are grateful to the organizers of the Summer School and
Conference on Poisson Geometry (2005) and  the International Center
of Theoretical Physics in Trieste of Italy, for their hospitality
while part of the work was being done.

\end{document}